
\input amssym.def
\input amssym.tex


\def\item#1{\vskip1.3pt\hang\textindent {\rm #1}}


\tolerance=300
\pretolerance=200
\hfuzz=1pt
\vfuzz=1pt


\hoffset=0.6in
\voffset=0.8in

\hsize=5.8 true in 


\vsize=8.5 true in
\parindent=25pt
\mathsurround=1pt
\parskip=1pt plus .25pt minus .25pt
\normallineskiplimit=.99pt

\countdef\revised=100
\mathchardef\emptyset="001F 
\chardef\ss="19
\def\3{\ss}
\def\anf{$\lower1.2ex\hbox{"}$}
\def\frac#1#2{{#1 \over #2}}
\def\>{>\!\!>}
\def\<{<\!\!<}

\def\ssarr{\hbox to 30pt{\rightarrowfill}}
\def\sarr{\hbox to 40pt{\rightarrowfill}}
\def\arr{\hbox to 60pt{\rightarrowfill}}
\def\larr{\hbox to 60pt{\leftarrowfill}}
\def\Arr{\hbox to 80pt{\rightarrowfill}}

{}

\def\ad{\mathop{\rm ad}\nolimits}

\def\Ad{\mathop{\rm Ad}\nolimits}

\def\conv{\mathop{\rm conv}\nolimits}



%
%

\def\Im{\mathop{\rm Im}\nolimits}



\def\Symm{\mathop{\rm Symm}\nolimits}
\def\Sp{\mathop{\rm Sp}\nolimits}



\def\0{{\bf 0}}
\def\1{{\bf 1}}

\def\a{{\frak a}}

\def\g{{\frak g}}

\def\k{{\frak k}}

\def\m{{\frak m}}

\def\n{{\frak n}}

\def\p{{\frak p}}

\def\z{{\frak z}}

\def\C{{\Bbb C}} 
\def\D{{\Bbb D}}

\def\N{{\Bbb N}}

\def\R{{\Bbb R}}

\def\:{\colon}  
\def\.{{\cdot}}
\def\|{\Vert}
\def\bsk{\bigskip}

\def\giantskip{\vskip2\bigskipamount}
\def\gsk{\giantskip}

\def\msk{\medskip}

\def\ssk{\smallskip}

\def\bbr{\bigbreak}
\def\giantbreak{\par \ifdim\lastskip<2\bigskipamount \removelastskip
         \penalty-400 \giantskip\fi}

\def\nin{\noindent}
\def\cen{\centerline}
\def\pagebreak{\vskip 0pt plus 0.0001fil\break}
\def\linebreak{\break}

\def\phi{\varphi}
\def\epsilon{\varepsilon}

\def\nin{\noindent}
\def\oline{\overline}

\def\pder#1,#2,#3 { {\partial #1 \over \partial #2}(#3)}
\def\pde#1,#2 { {\partial #1 \over \partial #2}}


\def\subeq{\subseteq}
\def\supeq{\supseteq}

\def\Rarrow{\Rightarrow}
\def\Larrow{\Leftarrow}

\def\up{{\uparrow}}

\font\eightrm=cmr8


\font\smc=cmcsc10
\font\bfone=cmbx10 scaled\magstep1 
\font\bftwo=cmbx10 scaled\magstep2 

\def\qed{{\unskip\nobreak\hfil\penalty50\hskip .001pt \hbox{}\nobreak\hfil
          \vrule height 1.2ex width 1.1ex depth -.1ex
           \parfillskip=0pt\finalhyphendemerits=0\medbreak}\rm}

\def\qeddis{\eqno{\vrule height 1.2ex width 1.1ex depth -.1ex} $$
                   \medbreak\rm}

\def\Lemma #1. {\bigbreak\vskip-\parskip\noindent{\bf Lemma #1.}\quad\it}

\def\Sublemma #1. {\bigbreak\vskip-\parskip\noindent{\bf Sublemma #1.}\quad\it}

\def\Proposition #1. {\bigbreak\vskip-\parskip\noindent{\bf Proposition #1.}
\quad\it}

\def\Corollary #1. {\bigbreak\vskip-\parskip\nin{\bf Corollary #1.}
\quad\it}

\def\Theorem #1. {\bigbreak\vskip-\parskip\noindent{\bf Theorem #1.}
\quad\it}

\def\Definition #1. {\rm\bigbreak\vskip-\parskip\noindent{\bf Definition #1.}
\quad}

\def\Remark #1. {\rm\bigbreak\vskip-\parskip\noindent{\bf Remark #1.}\quad}

\def\Example #1. {\rm\bigbreak\vskip-\parskip\noindent{\bf Example #1.}\quad}

\def\Problems #1. {\bigbreak\vskip-\parskip\noindent{\bf Problems #1.}\quad}
\def\Problem #1. {\bigbreak\vskip-\parskip\noindent{\bf Problems #1.}\quad}

\def\Conjecture #1. {\bigbreak\vskip-\parskip\noindent{\bf Conjecture #1.}\quad}

\def\Proof#1.{\rm\par\ifdim\lastskip<\bigskipamount\removelastskip\fi\smallskip
            \noindent {\bf Proof.}\quad}

\def\Axiom #1. {\bigbreak\vskip-\parskip\noindent{\bf Axiom #1.}\quad\it}

\def\Satz #1. {\bigbreak\vskip-\parskip\noindent{\bf Satz #1.}\quad\it}

\def\Korollar #1. {\bbr\vskip-\parskip\nin{\bf Korollar #1.} \quad\it}

\def\Bemerkung #1. {\rm\bigbreak\vskip-\parskip\noindent{\bf Bemerkung #1.}
\quad}

\def\Beispiel #1. {\rm\bigbreak\vskip-\parskip\noindent{\bf Beispiel #1.}\quad}
\def\Aufgabe #1. {\rm\bigbreak\vskip-\parskip\noindent{\bf Aufgabe #1.}\quad}

\def\Beweis#1. {\rm\par\ifdim\lastskip<\bigskipamount\removelastskip\fi
           \smallskip\noindent {\bf Beweis.}\quad}

\nopagenumbers

\def\date{\ifcase\month\or January\or February \or March\or April\or May
\or June\or July\or August\or September\or October\or November
\or December\fi\space\number\day, \number\year}

\def\title{Title ??}
\def\author{Author ??}

\def\thanks#1{\footnote*{\eightrm#1}}

\def\rightheadline{\hfil{\eightrm\title}\hfil\tenbf\folio}
\def\leftheadline{\tenbf\folio\hfil{\eightrm\author}\hfil}
\headline={\vbox{\line{\ifodd\pageno\rightheadline\else\leftheadline\fi}}}

\def\firstheadline{}
\def\firstfootline{\cen{\rm\folio}}

\def\seite #1 {\pageno #1
               \headline={\ifnum\pageno=#1 \firstheadline
               \else\ifodd\pageno\rightheadline\else\leftheadline\fi\fi}
               \footline={\ifnum\pageno=#1 \firstfootline\else{}\fi}}

\newdimen\dimenone
 \def\checkleftspace#1#2#3#4{
 \dimenone=\pagetotal
 \advance\dimenone by -\pageshrink   
 \ifdim\dimenone>\pagegoal          
   \else\dimenone=\pagetotal
        \advance\dimenone by \pagestretch
        \ifdim\dimenone<\pagegoal
          \dimenone=\pagetotal
          \advance\dimenone by#1         
          \setbox0=\vbox{#2\parskip=0pt                
                     \hyphenpenalty=10000
                     \rightskip=0pt plus 5em
                     \noindent#3 \vskip#4}    
        \advance\dimenone by\ht0
        \advance\dimenone by 3\baselineskip   
        \ifdim\dimenone>\pagegoal\vfill\eject\fi
          \else\eject\fi\fi}


\def\subheadline #1{\nin\bigbreak\vskip-\lastskip
      \checkleftspace{0.7cm}{\bf}{#1}{\medskipamount}
          \indent\vskip0.7cm\centerline{\bf #1}\medskip}

\def\sectionheadline #1{\bigbreak\vskip-\lastskip
      \checkleftspace{1.1cm}{\bf}{#1}{\bigskipamount}
         \vbox{\vskip1.1cm}\cen{\bfone #1}\bsk}

\def\lsectionheadline #1 #2{\bigbreak\vskip-\lastskip
      \checkleftspace{1.1cm}{\bf}{#1}{\bigskipamount}
         \vbox{\vskip1.1cm}\cen{\bfone #1}\msk \cen{\bfone #2}\bsk}

\def\lchapterheadline #1 #2{\bigbreak\vskip-\lastskip\indent\vskip3cm
                       \cen{\bftwo #1} \msk \cen{\bftwo #2} \gsk}
\def\llsectionheadline #1 #2 #3{\bigbreak\vskip-\lastskip\indent\vskip1.8cm
\cen{\bfone #1} \msk \cen{\bfone #2} \msk \cen{\bfone #3} \nobreak\bsk\nobreak}


\newtoks\literat
\def\[#1 #2\par{\literat={#2\unskip.}%
\hbox{\vtop{\hsize=.15\hsize\nin [#1]\hfill}
\vtop{\hsize=.82\hsize\nin\the\literat}}\par
\vskip.3\baselineskip}

\mathchardef\emptyset="001F 
\def\address{Author: \tt$\backslash$def$\backslash$address$\{$??$\}$}

\def\firstpage{\nin
{\obeylines \parindent 0pt }
\vskip2cm
\centerline {\bfone \title}
\gsk
\centerline{\bf\author}

\vskip1.5cm \rm}

\def\addresstwo{}

\def\dlastpage{\par\vbox{\vskip1cm\nin
\line{
\vtop{\hsize=.5\hsize{\parindent=0pt\baselineskip=10pt\nin\address}}
\quad 
\vtop{\hsize=.42\hsize\nin{\parindent=0pt
\baselineskip=10pt\addresstwo}}
\hfill} }}


\def\firstpage{\nin
{\obeylines \parindent 0pt }
\vskip2cm
\centerline {\bfone \title}
\ssk
\centerline {\bfone \titletwo}
\gsk
\centerline{\bf\author}
\vskip1.5cm \rm}

\def\bs{\backslash} 
\def\addots{\mathinner{\mkern1mu\raise1pt\vbox{\kern7pt\hbox{.}}\mkern2mu
\raise4pt\hbox{.}\mkern2mu\raise7pt\hbox{.}\mkern1mu}}

\pageno=1
\def\up#1{\leavevmode \raise.16ex\hbox{#1}}
 at 8truept
 at 8truept
 at 12truept
\chardef\ss="19
\def\3{\ss}
\def\title{Invariant Stein domains in Stein symmetric spaces}
\def\titletwo{and a non-linear complex convexity theorem}
\def\author{Simon Gindikin$^*$ and Bernhard Kr\"otz${}^\dagger$}
\footnote{}{${}^*$ Supported in part by the NSF-grant DMS-0070816 
and the MSRI}
\footnote{}{${}^\dagger$ Supported in part by the NSF-grant DMS-0097314 and the MSRI}

\def\date{October 11, 2001}
\def\Box #1 { \msk\par\nin 
\centerline{
\vbox{\offinterlineskip
\hrule
\hbox{\vrule\strut\hskip1ex\hfil{\smc#1}\hfill\hskip1ex}
\hrule}\vrule}\msk }

\def\address
{Simon Gindikin

Department of Mathematics

Rutgers University

New Brunswick, NJ 08903

USA

{\tt gindikin@math.rutgers.edu}

}

\def\addresstwo
{Bernhard Kr\"otz

The Ohio State University 

Department of Mathematics 

231 West 18th Avenue 

Columbus, OH 43210--1174 

USA

{\tt kroetz@math.ohio-state.edu}
}

\firstpage

\subheadline{Abstract}

We prove a complex version of Kostant's non-linear convexity theorem.
Applications to the construction of 
$G$-invariant Grauert tubes of Riemannian symmetric $G/K$ spaces are 
given.

\sectionheadline{Introduction}

Let $X=G/K$ be a semisimple non-compact Riemannian symmetric space.  
We may assume that $G$ is semisimple with finite center and we write $G=NAK$ 
for an Iwasawa decomposition of $G$. By our assumption,  $G$ 
sits in its universal complexification $G_\C$ and so 
$X\subeq X_\C\:=G_\C/K_\C$. Note that $X_\C$ is a Stein symmetric 
space. Observe that the group $G$ does not act properly 
(i.e. with compact isotropy subgroups) on $G_\C/K_\C$ since 
$K_\C$ is not compact.  Write $U$ for the compact real 
form of $G_\C$. The maximal
connected subdomain of $GUK_\C/ K_\C\subeq X_\C$ on which $G$ acts properly
and which contains $X$ was introduced in [AkGi90]. It is given by 
$$\Xi\:=G\exp(i\Omega)K_\C/ K_\C$$
where $\Omega$ is a polyhedral convex domain in $\a\:={\rm Lie}(A)$
defined by 
$$\Omega\:=\{X\in\a\: (\forall\alpha\in\Sigma)\ |\alpha(X)|<{\pi\over 2}\}.$$
Here $\Sigma$ denotes the restricted root system with respect to $\a$.

\par One of our principal aims is to construct a broad class 
of $G$-invariant Stein subdomains in $\Xi$. Let us remind the reader 
that it is a notoriously  difficult problem to verify that 
certain unions of orbits of non-compact groups are Stein. 
For example the problem posed in [AkGi90] whether $\Xi$ is Stein was unsolved 
until the last year. 

\par The Iwasawa decomposition on $G$ cannot be holomorphically 
extended to the whole group $G_\C$; we have that $N_\C A_\C K_\C\subsetneq G_\C$ is 
an open Zariski dense subset. The {\it Iwasawa domain}
$\Xi_I$ is defined as  the maximal connected $G$-invariant subdomain
in $X_\C$ which contains $X$ and is contained in $N_\C A_\C K_\C/K_\C$:
$$\Xi_I\:=\big(\bigcap_{g\in G} g(N_\C A_\C K_\C/ K_\C)\big)_0,$$
where $(\cdot)_0$ refers to the connected component containing 
$X$. Since $\Xi_I$ is the connected component 
of an open intersection of Stein domains, it is easy 
to see that $\Xi_I$ is Stein (cf.\ proof of Theorem 3.4).

\par The domains $\Xi$ and $\Xi_I$ again became of recent interest. 
One has $\Xi=\Xi_I$. Here are the results in chronological order:

\msk 
\item{${\bf\cdot}$} $\Xi\subeq \Xi_I$ for all classical groups $G$
(cf.\ [KrSt01a]).
\item{${\bf\cdot}$} $\Xi=\Xi_I$ for all classical groups $G$
(cf.\ [GiMa01]).
\item{${\bf\cdot}$} $\Xi_I\subeq \Xi$ for all $G$
(cf.\ [Ba01]).

\msk 
In [KOS01] it will be shown that $\Xi\subeq \Xi_I$ for the 
exceptional cases. Also it is 
announced in [Hu01] that $\Xi\subeq \Xi_I$. This will
then give us $\Xi=\Xi_I$ in the general case. In particular, 
$\Xi$ is Stein. Also, in the  preprint 
[BHH01] a complex-geometric proof of the Steinness of $\Xi$ is given. 
 
\ssk Geometrically, the domain $\Xi$ is a rather complicated object. 
If $G$ is a group of Hermitian type, then $\Xi=X\times \oline X$
(cf.\ [BHH01], [GiMa01] or [KrSt01b]). 
If $X=G/K$ is classical, then there exists a 
group of Hermitian type $S\supeq G$ with maximal compact subgroup 
$U\supeq K$ such that 
$G/K$ is a real form of the Hermitian symmetric space $S/U$. 
In [BHH01] and [KrSt01b] a $G$-equivariant subdomain $\Xi_0\subeq \Xi$
was exhibited which is biholomorphic to $S/U$. Further one has
$\Xi=\Xi_0=S/U$ if and only if $\Sigma$ is of type $C_n$ or $BC_n$. 
In particular, if $\Sigma\neq C_n, BC_n$, then 
$S/U\subsetneq \Xi$ and the explicit geometric structure  of $\Xi$ 
is very intricate. For many exceptional spaces $G/K$ the geometric 
structure is even more complicated.  

\ssk The domain $\Xi$ (in [Gi98] it is called {\it 
complex crown of} $X$) is universal in the sense that many 
analytical and geometrical constructions on the Riemannian symmetric 
space $X$ extend to $\Xi$. 
In [KrSt01a] it was shown that
the inclusion $\Xi\subeq \Xi_I$ implies that  
all eigenfunctions on $X$ for the algebra of $G$-invariant 
differential operators $\D(X)$ extend holomorphically to $\Xi$. 
In [GiMa01] the problem $\Xi\subeq \Xi_I$ is included in a broad 
class of geometrical problems connected with Matsuki duality. 
These geometrical problems include in particular the problem of 
the parametrization of compact complex cycles in flag domains
(cf.\ [Wo92]).

\msk Let us now come to the contents of this paper. Our first main 
result is a complex version of Kostant's non-linear convexity theorem 
(cf.\ [Kos73]). Write $G\to A,\  g\mapsto a(g)$ for the 
middle projection in $G=NAK$. 

\msk\nin {\bf Theorem.} (Kostant)  {\it Let $X\in \a$. Then 
$$a(K\exp(X))=\exp(\conv({\cal W}X)),$$
where ${\cal W}=N_K(\a)/Z_K(\a)$ denotes the Weyl group 
and $\conv(\cdot)$ refers to the convex hull of $(\cdot)$.}
\msk 

One shows that the  middle projection $a\: G\to A$ holomorphically 
extends to 
$$a\: G\exp(i\Omega)K_\C \to A_\C.\leqno(1)$$
Then our main result is:

\msk 
\nin {\bf Theorem A.} (Complex Convexity Theorem)
{\it For all $X\in \Omega$ we have that 
$$a(G\exp(iX))\subeq A\exp(i\conv({\cal W}X)).$$}
\msk

Our convexity theorem features interesting applications to  the 
geometry of the domain $\Xi$ and its generic 
subdomains which are defined as follows: Let $\omega\subeq \Omega$
be a non-empty convex ${\cal W}$-invariant open subset. Then we can form 
the domains

$$\Xi(\omega)\:=G\exp(i\omega)K_\C/ K_\C.$$
Note that $\Xi=\Xi(\Omega)$ and that the Iwasawa projection 
(1) naturally factors to a holomorphic mapping $a\: \Xi\to A_\C$. 
Finally we define for every $g\in G$ 
the {\it horospherical tube} $T(g,\omega)\subeq X_\C$
by 
$$T(g,\omega)\:= g\big(N_\C A\exp(i\omega)K_\C/K_\C\big).$$ 

\par As an application of Theorem A we now obtain: 

\msk\nin {\bf Theorem B.} {\it Let $\omega$ be an open convex 
Weyl group invariant subset of $\Omega$. Then the following assertions hold:
\item{(i)}  $a(\Xi(\omega))=A\exp(i\omega)$.
\item{(ii)} $\Xi(\omega)=\big(\bigcap_{g\in G} T(g,\omega)\big)_0$. 
\item{(iii)} The domain $\Xi(\omega)$ is Stein.}

\msk Let us emphasize that in the the case of $\Xi(\omega)=\Xi$, 
the inclusion $\Xi\subeq \Xi_I$ means only the existence 
of an Iwasawa projection $a\: \Xi\to A_\C$; (i) in Theorem B  gives
a more precise information on the image of this projection. 
Further (ii) in the above theorem is a much stronger statement than 
(iii); in particular, (ii) implies (iii) since 
all horospherical tubes $T(g,\omega)$ are Stein. 
 
\msk Theorem B can be considered as an analogue 
of Lassalle's results for compact symmetric spaces (cf.\ [La78]). We can 
interpret Theorem A and Theorem B as statements for $G$-orbits 
in $X_\C$ intersecting  $A\exp(i\Omega)K_\C/ K_\C$ which in 
in the case of  compact symmetric spaces $G/K$ are 
true for arbitrary $G$-orbits (cf.\ [La78]).

\msk It is our pleasure to thank the MSRI, Berkeley,  for 
its hospitality during the {\it Integral geometry program}
where this work was accomplished. We  thank 
Dmitri Akhiezer for his careful screening of the manuscript 
and his worthy suggestions. Further we would like to
thank Karl-Hermann Neeb for his kindness to proofread the paper  
and Laura Geatti for sharing her knowledge with us on the boundary 
of $\Xi$.

\sectionheadline{1. Notation} 

Let $G$ be a connected semisimple Lie group sitting inside a 
complexification $G_\C$. 
We denote by $\g$ and $\g_\C$ the Lie algebras of $G$ and $G_\C$, 
respectively. 
Let $K<G$ be a maximal compact subgroup and $\k$ its Lie algebra. Denote by 
$\theta\: G\to G$ a Cartan involution which has $K$ as a fixed point set. 

\par   Let $\g=\k\oplus\p$ the Cartan decomposition attached to $\k$. 
Take  $\a\subeq 
\p$ a maximal Abelian subspace and let $\Sigma=\Sigma(\g, \a)\subeq \a^*$ 
be the 
corresponding root system. Related to this root system is the root space 
decomposition according to the simultaneous eigenvalues of $\ad (H), H\in \a:$

$$\g=\a\oplus \m \oplus\bigoplus_{\alpha\in \Sigma} \g^\alpha,$$ 
here $\m=\z_\k(\a)$ and $\g^\alpha =\{ X\in \g\: (\forall H\in \a) \ 
[H, X]=\alpha(H) X\}$. For the choice of a positive system $\Sigma^+\subeq 
\Sigma$ 
one obtains the nilpotent Lie algebra $\n =\bigoplus_{\alpha\in 
\Sigma^+}\g^\alpha$. Then one has the Iwasawa decomposition on the Lie algebra level 
$$\g=\n\oplus \a \oplus \k.$$ 
We write $A$, $N$ for the analytic subgroups of $G$ corresponding to 
$\a$ and $\n$.  
For these choices one has for $G$ the Iwasawa decomposition, namely, 
the multiplication map 

$$N\times A\times K\to G, \ \ (n, a, k)\mapsto nak$$
In particular, every element 
$g\in G$ can be written uniquely as $g=n(g) a(g) \kappa(g)$ with each of the 
maps $\kappa(g)\in K$, $a(g)\in A$, $n(g)\in N$ depending analytically on $g\in 
G$. 
The last piece of structure theory we shall recall is the little Weyl group. 
We 
denote by ${\cal W} =N_K(\a)/ Z_K(\a)$ the {\it Weyl group} of 
$\Sigma(\a,\g)$. 

Finally we define the domain
$$\Omega=\{ X\in\a\:(\forall \alpha\in \Sigma)\  |\alpha(X)|<{\pi\over 2}\}.$$
Clearly $\Omega$ is convex and ${\cal W}$-invariant.

\sectionheadline{2.  The complex convexity theorem}

\ssk  Let us first give the relevant notation. 
Denote by $K_\C$, $A_\C$, $N_\C$
the complexifications of $K$, $A$ and $N$. 
Then $N_\C A_\C K_\C$ is a proper Zariski-open, hence  dense subset 
of $G_\C$. 
Throughout this paper we will assume that $\Xi\subeq \Xi_I$ holds, i.e.,  

$$ G\exp(i\Omega) \subeq N_\C A_\C K_\C$$
(see the introduction for proofs of this statement).  
We now set $T_\Omega\:=A\exp(i\Omega)\subeq A_\C$. 
Then $G T_\Omega \subeq  N_\C A_\C K_\C$.

\par One can show that one has a well defined holomorphic 
middle projection $N_\C A_\C K_\C \to A_\C/ (A_\C\cap K_\C)$. 
Due to the simple connectedness of $T_\Omega$, this projection 
restricted to $G T_\Omega $ lifts to $A_\C$ (cf.\ [KrSt01a, proof of 
Th.\ 1.8(iii)]) and we obtain an analytic  mapping 

$$G\times T_\Omega \to A_\C, \ \ (g,a)\mapsto a(ga)$$
holomorphic in the second variable such that 

$$ga\in N_\C a(ga) K_\C$$
holds. 
  
\par If $V$ is a vector space and $E\subeq V$ is a subset, then we 
denote by $\conv(E)$ the convex hull of $E$ in $V$.

\Theorem 2.1. {\rm (Complex Convexity Theorem)} 
Assume that $G$ is classical.  Then  we have for all $X\in \Omega$ that 

$$a(G\exp(iX))\subeq A\exp(i\conv({\cal W}X)).\qeddis

The proof of Theorem 2.1 will be given in several steps.

Fix $a\in T_\Omega$ and consider the function

$$f_a\: K\to \a_\C, \ \ k\mapsto \log a(ka).$$ 
Further we write $p_{\a_\C}\: \g_\C \to \a_\C$ for the 
projection along $\k_\C +\n_\C$. For  $x\in GT_\Omega K_\C \subeq 
N_\C A_\C K_\C$ we write $b(x)=n(x) a(x)$ for the {\it triangular  part}
of $x$.

\Lemma 2.2. For any $a\in T_\Omega$, $k\in K$ and $X\in\k$ we have 

$${d\over dt}\Big|_{t=0} f_a(\exp(tX)k)= p_{\a_\C} (\Ad(b(ka))^{-1}X).$$

\Proof. This result can be easily deduced from the known case for $a\in A$
by analytic continuation. However, for conveniene for the reader, 
we briefly recall the proof. 
We have 

$$\eqalign{{d\over dt}\Big|_{t=0} f_a(\exp(tX)k)
&={d\over dt}\Big|_{t=0}\log a(\exp(tX)ka)\cr 
&={d\over dt}\Big|_{t=0}\log a(\exp(tX)b(ka))\cr 
&= 
{d\over dt}\Big|_{t=0}[\log a(b(ka)^{-1}\exp(tX)b(ka))+\log a(ka)]\cr 
&={d\over dt}\Big|_{t=0}\log 
a(\exp (t\Ad (b(ka))^{-1}X))= p_{\a_\C} (\Ad(b(ka))^{-1}X),\cr}$$
proving the lemma.\qed

Write $\g_\C^\R$ for $\g_\C$ considered as a real Lie algebra and $\kappa_\R$
for the Cartan-Killing form on  $\g_\C^\R$. Let $\kappa$ be  the 
Cartan-Killing form on $\g$ and recall the following relation 
between $\kappa$ and $\kappa_\R$:

$$(\forall X,X',Y,Y'\in\g)\qquad \kappa_\R(X+iY, X'+iY')=
2(\kappa(X,X')-\kappa(Y,Y')).$$ 
For every $\lambda\in (\g_\C^\R)^*$ we define $H_\lambda\in \g_\C^\R$
by $\lambda(X)=\kappa_\R(X,H_\lambda)$ for all $X\in\g_\C^\R$. 

For every $\lambda\in (\a_\C^\R)^*$ and $a\in T_\Omega$ we now define 
the function 

$$f_{a,\lambda}\: K\to\R, \ \ k\mapsto \lambda(f_a(k)).$$

\Lemma 2.3. Let $a\in T_\Omega$ and $k\in K$. 
Then the following assertions hold: 
\item{(i)} For all $X\in\k$ one has 
${d\over dt}\Big|_{t=0} f_{a,\lambda}(\exp(tX)k)= 
\kappa_\R (X, \Ad(n(ka))H_\lambda).$
\item{(ii)} We have $df_{a,\lambda}(k)=0$ if and only if 
$$\k\bot_{\kappa_\R} \Ad(n(ka))H_\lambda.$$

\Proof. (i) is immediate from Lemma 2.2 and the notations
introduced from above. 
\par\nin (ii) follows from (i).\qed 

Write $i \a_\R^*$ for the subspace of $(\a_\C^\R)^*$ which vanishes 
on $\a$. Our next goal is to determine the critical set of $f_{a,\lambda}$ 
for $\lambda\in i\a_\R^*$. 

\par We also denote by $\theta$ the holomorphic extension of the 
Cartan involution to $G_\C$. Further we write $G_\C\to G_\C,
\ x\mapsto \oline x$ for the complex conjugation with respect 
to the real form $G$. 

\Lemma 2.4. Suppose that $X\in \Omega$ is regular. Set $a=\exp(iX)$
and let $k\in K$. Then 
the following implication holds
$$ka\in N A_\C K_\C \Rarrow k\in N_K(\a).$$

\Proof. We write $ka=n bk'$ for $n\in N$, $b\in A_\C$ and 
$k'\in K_\C$. Then 

$$ka^2k^{-1}=ka(\theta(ka)^{-1})=n bk'(\theta(nbk'))^{-1}=
nb^2\theta(n)^{-1}.$$
Set $x\:=ka^2k^{-1}$. Then 
$$\oline x= ka^{-2}k^{-1}=x^{-1}.$$
Thus we obtain that 
$$\oline x= n\oline b^2 \theta(n)^{-1}= \theta(n) b^{-2} n^{-1}.$$
Therefore $x= \theta(n) \oline b^{-2} n^{-1}$ and so 
$$x^2=\theta(n) \oline b^{-2} n^{-1} n b^2\theta(n)^{-1}\in A_\C 
\theta(N_\C).$$
Now let $0<t\leq 1$ and write $a_t\:=\exp(itX)$ and accordingly 
we define $x_t$. Then we obtain that 

$$\log x_t^2=\log(k\exp(i4tX)k^{-1})=\Ad(k)(i4tX)\in \a_\C +\theta(\n_\C)$$
for small $t$ and hence for all $t$ by analyticity. But 
$\Ad(k)(i4tX)\in i\p$. Hence  by the regularity of $X$ we obtain 
that $k\in N_K(\a)$, as was to be shown. \qed

\Proposition 2.5.  Suppose that $X\in \Omega$  and set $a=\exp(iX)$. 
Let $\lambda\in i\a_\R^*$ be such that $H_\lambda\in i\a$ is regular. 
Then we have 

$$df_{a,\lambda}(k)=0 \iff k\in N_K(\a).$$

\Proof. The implication ``$\Larrow$'' follows immediately 
from Lemma 2.3 (ii). 
\par\nin ``$\Rarrow$'': Suppose that $df_{a,\lambda}(k)=0$. Then 
Lemma 2.3 (ii) implies that $\k\bot_{\kappa_\R} \Ad (n(ka))H_\lambda$, 
or equivalently $ \Ad (n(ka))H_\lambda\in i\g+ \p$.  
Assume that $k\not\in N_K(\a)$. By Lemma 2.4 we have 
$n(ka)=\exp(Y)$ for some $Y\in \n_\C\bs \n$. Hence 
$$\Ad(n(ka))H_\lambda= H_\lambda +\underbrace{[Y,H_\lambda]}_{\in\n_\C}
+{1\over 2} \underbrace{[Y,[Y,H_\lambda]]}_{\in [\n_\C,\n_\C]}+\ldots.$$
Thus the fact that $H_\lambda\in i\a$ is regular
implies that $\Ad(n(ka))H_\lambda\in i\a +(\n_\C\bs i\n)$. But 
this contradicts $ \Ad (n(ka))H_\lambda\in i\g+ \p$, 
proving our proposition.\qed

\par\nin {\bf Proof of Theorem 2.1.} First we observe that it is sufficient 
to prove 
$$a(K\exp(iX))\subeq A\exp(i\conv({\cal W}X))$$
for all $X\in \Omega$. By a simple density/continuity  argument 
we  may further assume that $X$ is regular.
Suppose that there exists a $k\in K$ such that 
$a(k\exp(iX))\not\in  A\exp(i\conv({\cal W}X)),$ or equivalently 
$$\Im \log a(k\exp(iX))\not \in \conv({\cal W}X).$$
Then we find a regular element $\lambda\in i\a_\R^*$ such that 
$$f_{\exp(iX),\lambda}(k)> \max_{Y\in i \conv({\cal W}X)} \lambda(Y).$$
But Proposition 2.5 implies that $f_{\exp(iX),\lambda}$ takes its 
maximum at an element $k\in N_K(\a)$. Hence 
$f_{\exp(iX),\lambda}(k)=\lambda(\Ad(k)iX)$; a contradiction 
to our inequality from above.\qed

\sectionheadline{3. Applications}

Let now $\omega\subeq \Omega$ be an open Weyl group invariant 
convex set. Then we define the domain 
$$\Xi(\omega)\:= G\exp(i\omega) K_\C/ K_\C. $$
If $\omega=\Omega$, then we set $\Xi\:=\Xi(\Omega)$. 
Write $\partial \Xi(\omega)$ for the toplogical 
boundary of $\Xi(\omega)$ in $G_\C/K_\C$. 
Note the following properties of $\Xi(\omega)$: 

\msk 
\item{${\bf\cdot}$} $\Xi(\omega)$ is open in $G_\C/ K_\C$ (cf.\ [AkGi90]). 
\item{${\bf\cdot}$} $\Xi(\omega)$ is connected and $G$-invariant. 
\item{${\bf\cdot}$} $G$ acts properly on $\Xi$ (cf.\ [AkGi90]).  
\item{${\bf\cdot}$} One has $G\exp(i\partial \omega)K_\C/ K_\C
\subeq \partial \Xi(\omega)$. Moreover if $\oline\omega\subeq 
\Omega$, then we have $\partial \Xi(\omega)=G\exp(i\partial \omega)
K_\C/ K_\C$ (cf.\ [KrSt01b, Prop.\ 4.1]). 

\msk\nin  From our discussions in the previous section it is clear 
that we have a holomorphic projection 
$$a\: \Xi_\omega\to A_\C,\ \  x\mapsto a(x)$$
with $x\in N_\C a(x)K_\C/ K_\C$ for all $x\in \Xi$.   
Finally we define the Abelian tube domain 
$$T_\omega\:=A\exp(i\omega)\subeq A_\C.$$

An immediate consequence of our complex convexity theorem then is: 

\Lemma 3.1.  We have 
$$a(\Xi(\omega))=T_\omega.\qeddis 

\Remark 3.2. From Lemma 3.1 we obtain in particular that 
$a(\Xi)\subeq T_\Omega$. It is interesting to observe that this 
inclusion for $G=\Sp(n,\R)$ extends a result of Siegel. 
\par Consider the vector space $V\:=\Symm(n,\R)$ of real symmetric 
matrices  with its 
subcone of positive definite matrices $V_+$. Then we have 
the symmetric Siegel domain
$$S^+\:=V+iV_+\subeq V_\C.$$
Recall that $S^+\cong \Sp(n,\R)/U(n)$. 
\par If we write $\Delta_j$ for the $j$-th principal minor on $V_\C$, 
then Siegel's Lemma says
$$\Delta_j(z)\neq 0 \qquad\hbox{for $z\in S^+$}\leqno(3.1)$$ 
and all $1\leq j\leq n$. Now consider the rational functions
$$\chi_j(z)={\Delta_j(z)\over \Delta_{j-1}(z)}$$
on $V_\C$. Then the inclusion $a(\Xi)\subeq T_\Omega$
for $G=\Sp(n,\R)$ implies
$$\Im \chi_j(z)>0 \qquad \hbox{for $z\in S^+$}\leqno(3.2)$$
and for all $j$. Note that (3.2) implies (3.1). To see this 
one identifies $S^+$ with the symmetric space $\Sp(n,\R)/U(n)$. 
Then with $S^-\:=\oline {S^+}$
one has a biholomorphism $\Xi\cong S^+\times S^-$ (cf.\ [GiMa01]
or [KrSt01b]). Realizing $A$ as diagonal matrices in $G$, 
one then easily shows that $a(\Xi)\subeq T_\Omega$
implies (3.2).
\par The example discussed above admits a natural generalization 
to all   tube domains $V+iV_+$ associated to an Euclidean Jordan algebra
$V$ and cone $V_+$. \qed

For every $g\in G$ we define the 
{\it horospherical tube} associated to $\omega$ by 

$$T(g,\omega)\:= g\big( N_\C T_\omega K_\C/K_\C\big)\subeq G_\C/K_\C.$$
Note that $T(\1,\omega)$ is biholomorphic to $N_\C\times T_\omega$. 
In particular, all horospherical tubes $T(g,\omega)$ are Stein. 

\Theorem 3.3. For any  non-empty 
open convex ${\cal W}$-invariant set $\omega\subeq \Omega$ the domain 
$\Xi(\omega)$ is the connected component of the 
intersection of horospherical tubes:
$$\Xi(\omega)=\big(\bigcap_{g\in G}T(g,\omega)\big)_0.$$

\Proof. From Lemma  3.1 we obtain that 
$$\Xi(\omega)\subeq N_\C T_\omega K_\C/ K_\C=a^{-1}(T(\1,\omega)).$$
Thus the fact that $\Xi(\omega)$ is $G$-invariant 
and connected implies that 
$\Xi(\omega)\subeq \big(\bigcap_{g\in G}T(g,\omega)\big)_0.$
On the other hand we obtain from $\Xi=\Xi_I$ (cf.\ the discussion 
in the introduction) and Lemma 3.1 
that $\Xi=\big(\bigcap_{g\in G}T(g,\Omega)\big)_0$.
{}From this and Lemma 3.1 we hence get 
$$\big(\bigcap_{g\in G}T(g,\omega)\big)_0\subeq \Xi\cap T(\1,\omega)
=\Xi(\omega),$$
completing the proof of the theorem.\qed 

An interesting  application of Theorem 3.3 is the following:

\Proposition 3.4. Let $\omega\subeq\Omega$ be an 
open convex ${\cal W}$-invariant set. Then the intersection 
$$I(\omega)\:=\bigcap_{g\in G} T(g,\omega)$$
is open. In particular, $I(\omega)$ and every connected component
of $I(\omega)$ is Stein. In particular, $\Xi(\omega)$ is Stein. 

\Proof. (following a suggestion of Dmitri Akhiezer) 
Since $G=KAN$ and since $AN$ leaves $T(\1,\omega)$ invariant, 
we obtain 
$$I(\omega)=\bigcap_{k\in K} T(k,\omega).$$
Hence $I(\omega)$ is an intersection of open sets over the ``compact parameter
space'' $K$. In particular $I$ is open. Since all horospherical 
tubes $T(k,\omega)$ are Stein, we obtain that $I$ is Stein. 
With $I(\omega)$ all its connected components are Stein, concluding the 
proof of the theorem. \qed 

As a final application of the complex convexity theorem 
we prove a result on the characterization of the boundary
of $\Xi$. 

\Proposition 3.5. Let $\omega\subeq\Omega$ be an 
open convex ${\cal W}$-invariant set. Let $(z_n)_{n\in\N}$ be a sequence with 
$z_n\to z_0\in \partial\Xi(\omega)$. Then
$$(\Im \log a(z_n))_{n\in \N}$$
is a sequence in $\omega$ and every accumulation point of this 
sequence lies in $\partial\omega$.

\Proof. Set $X_n\:=\Im \log a(z_n)$. By Lemma 3.1 we have $X_n\in \omega$.
Since $\oline \omega$ is compact in 
$\a$, we may assume that $X_n\to X_0$ with $X_0\in \oline\omega$.
\par If $X_0\not\in \partial\omega$, then we find a convex Weyl group
invariant open set $\omega_1$ such that $\oline{\omega_1}\subeq 
\omega$ and $X_n\in \omega_1$.
Then Theorem 3.3 implies that $z_n\in \Xi(\omega_1)$ and 
so $z_0\in \oline {\Xi(\omega_1)}$. Now $\oline \omega_1\subeq\Omega$ and 
so $\partial\Xi(\omega_1)=G\exp(i\partial\omega_1)K_\C/K_\C$. 
Thus $z_0\in  \Xi(\omega)$, contradicting the assumption 
$z_0\in \partial\Xi(\omega)$. This concludes the  proof of 
the proposition.\qed

\def\entries{

\[AkGi90 Akhiezer, D.\ N., and S.\ G.\ Gindikin, {\it On Stein extensions of 
real symmetric spaces}, 
Math.\ Ann.\ {\bf 286}, 1--12, 1990

\[Ba01 Barchini, L., {\it Stein Extensions of Real Symmetric Spaces 
and the Geometry of the Flag Manifold}, preprint

\[BHH01 Burns, D., S.\ Halverscheid, and R.\ Hind, {\it The Geometry 
of Grauert Tubes and Complexification of Symmetric Spaces}, 
preprint

\[Gi98 Gindikin, S., {\it Tube domains in Stein symmetric spaces}, 
Positivity in Lie theory: open problems, 81--97, de Gruyter Exp. Math., 
{\bf 26}, de Gruyter, Berlin, 1998

\[GiMa01 Gindikin, S., and T. Matsuki, {\it Stein Extensions
of Riemann Symmetric Spaces and Dualities of Orbits on Flag 
Manifolds}, MSRI preprint 2001-028

\[Hu01 Huckleberry, A., {\it On certain domains in cycle spaces
of flag manifolds}, preprint

\[Kos73 Kostant, B., {\it On Convexity, the Weyl Group and the Iwasawa 
Decomposition}, 
Ann.\ scient.\ \'Ec.\ Norm.\ Sup.\ $4^e$ s\'erie, t.\ {\bf 6}, 413--455, 1973 

\[KrSt01a Kr\"otz, B., and R.\ J. Stanton, 
{\it Holomorphic extensions of representations: (I)  
automorphic functions}, preprint

\[KrSt01b ---, {\it Holomorphic extensions of representations: (II)  
geometry and harmonic analysis}, preprint

\[KOS01 Kr\"otz, B., M. Otto, and R.\ J. Stanton, {\it Complex crowns of 
Riemannian symmetric spaces -- excercises for the exceptional cases}, 
in preparation 

\[La78 Lassalle, M., {\it S\'eries de Laurent des fonctions holomorphes 
dans la complexification d'un espace sym\'etrique compact}, 
Ann. Sci. \'Ecole Norm. Sup. {\bf (4) 11} (1978), no. {\bf 2}, 167--210 

\[Wo92 Wolf, J.\ A., {\it  The Stein condition for cycle spaces of 
open orbits on complex flag manifolds}, Ann. of Math. {\bf (2) 136} (1992),
no. {\bf 3}, 541--555

}

{\sectionheadline{\bf References}
\frenchspacing
\entries\par}
\dlastpage 
\bye
\end